\documentclass{article}
\usepackage{amsmath}
\usepackage{amsfonts}
\usepackage{amssymb}
\usepackage{graphicx}
\usepackage{cite}

\newtheorem{theorem}{Theorem}[section]

\newtheorem{definition}[theorem]{Definition}
\newtheorem{example}[theorem]{Example}

\newtheorem{lemma}[theorem]{Lemma}

\input{tcilatex}
\begin{document}

\title{On the Perturbation-Iteration Algorithm for System of Fractional Differential Equations}
\author{Mehmet \c{S}ENOL and \.{I}. Timu\c{c}in DOLAPC\.{I} \\
Nev\c{s}ehir Hac\i\ Bekta\c{s} Veli University, Department of Mathematics,
Nev\c{s}ehir, Turkey\\
Celal Bayar University, Department of Mechanical Engineering,\\
Manisa, Turkey\\
e-mail:msenol@nevsehir.edu.tr, ihsan.dolapci@cbu.edu.tr}
\maketitle

\begin{abstract}
In this study, perturbation-iteration algorithm, namely PIA, is applied to
solve some types of system of fractional differential equations (FDEs) for
the first time. To illustrate the efficiency of the method, numerical
solutions are compared with the results published in the literature by
considering some FDEs. The results confirm that the PIA is powerful, simple
and reliable method for solving system of nonlinear fractional differential
equations.

\textbf{Keywords:} Fractional-integro differential equations, Caputo
fractional derivative, Initial value problems, Perturbation-Iteration
Algorithm.
AMS 2010: 34A08, 34A12.
\end{abstract}

\section{Introduction}

\label{intro} As an important mathematical branch investigating the
properties of derivatives and integrals, the history of the fractional
calculus is nearly as old as classical integer order analysis. Even though
there is a long period of time since the beginning of fractional calculus,
it could not find practical applications for many years. However, in last
decades, it has found a place in various areas such as control theory, Senol
et al. (2014), viscoelasticity, Yu and Lin. (1998), electrochemistry, Oldham
(2010) and electromagnetic, Heaviside (2008).

The evolution of the symbolic computation programs such as Matlab and
Mathematica is one of the driving forces behind this increased usage.
Various multidisciplinary problems could be expressed by the help of
fractional derivatives and integrals. For example, in Hamamci, (2007). The
most important descriptions of fundamentals of fractional calculus have been
studied by Mainardi, (1997) and Podlubny, (1998). Existence and uniqueness
of the solutions has also been studied by Yakar and Koksal (2012) and the
references therein.

Parallel to the studies in applied sciences, system of fractional
differential equations (FDEs) allowed scientists to describe and model
various important and useful physical problems.

The number of differential equations whose solution can not be found
analytically. That situations appear in FDEs more than the other types of
differential equations. In this case, as the study of algorithms using
numerical approximation for the problems of mathematical analysis, the field
of numerical analysis is used for approximate solutions of FDEs. In recent
years, a significant effort has been expended to propose numerical methods
for this purpose. These methods include, fractional variational iteration
method, Wu and Lee, (2010); Guo and Mei (2011), homotopy perturbation
method, Abdulaziz et al. (2008); He (2012); Momani and Odibat, (2007); Zhang
et al. (2014) and fractional differential transform method (Momani et al.,
2007; Arikoglu and Ozkol, 2009; El-Sayed et al., 2014).

In this study, we have applied the previously developed method PIA to obtain
approximate solutions for some system of FDEs. Our method is suitable for a
broad class of equations and does not require special assumptions and
restrictions.

In the literature, there exists a few fractional derivative definitions of
an arbitary order. Two mostly used of them are the Riemann-Liouville and
Caputo fractional derivatives. The two definitions are quite similar but
have different order of evaluation of derivation.

\bigskip Riemann-Liouville fractional integration of order $\alpha $ is
defined by:%
\begin{equation}
J^{\alpha }f(x)=\frac{1}{\Gamma (\alpha )}\int_{0}^{x}(x-t)^{\alpha
-1}f(t)dt,\quad \alpha >0,\quad x>0.  \label{1}
\end{equation}

\qquad The following two equations are defined as Riemann-Liouville and
Caputo fractional derivatives of order $\alpha ,$ respectively.%
\begin{equation}
D^{\alpha }f(x)=\frac{d^{m}}{dx^{m}}\left( J^{m-\alpha }f(x)\right)
\label{2}
\end{equation}%
\begin{equation}
D_{\ast }^{\alpha }f(x)=J^{m-\alpha }\left( \frac{d^{m}}{dx^{m}}f(x)\right) .
\label{3}
\end{equation}%
where $m-1< \alpha < m $ and $m \in \mathbf{N}.$

Due to the appropriateness of the initial conditions, fractional definition
of Caputo is often used in recent years. \vskip .1in 
\begin{definition}
The fractional derivative of $u(x)$ $in$ the Caputo sense is defined as%
\begin{equation}
D_{\ast }^{\alpha }u(x)=\left\{ 
\begin{array}{cc}
\frac{1}{\Gamma (m-\alpha )}\int_{0}^{x}(x-t)^{m-\alpha -1}u^{(m)}(t)dt & 
,m-1<\alpha <m \\ 
\frac{d^{m}}{dx^{m}} &, \alpha =m%
\end{array}%
\right.  \label{4}
\end{equation}%
for $m-1<\alpha <m,$ $m\in  \mathbf{N},$ $x>0,$ $u\in C_{-1}^{m}.$
\end{definition}

Also, we need here two of its basic properties. \vskip .1in 
\begin{lemma}
If $m-1<\alpha <m,$ $m\in \mathbf{N}$ and $u\in C_{\mu }^{m},$ $\mu -1$ then 
\begin{equation}
D_{\ast }^{\alpha }J^{\alpha }u(x)=u(x)  \label{5}
\end{equation}
and 
\begin{equation}
J^{\alpha }D_{\ast }^{\alpha }u(x)=u(x)-\sum_{k=0}^{m-1}u^{(k)}(0^{+})\frac{%
x^{k}}{k!},\quad x>0.  \label{6}
\end{equation}
\end{lemma}

After this introductory section, Section 2 is reserved to a brief review of
the Perturbation-Iteration Algorithm PIA(1,1), in Section 3 some examples
are presented to show the efficiency and simplicity of the algorithm.
Finally the paper ends with a conclusion in Section 4.

\section{Overview of the Perturbation-Iteration Algorithm PIA(1,1)}

As one of the most practical subjects of physics and mathematics,
differential equations create models for a number of problems in science and
engineering to give an explanation for a better understanding of the events.
Perturbation methods have been used for this purpose for over a century
(Nayfeh, 2008; Jordan and Smith, 1987; Skorokhod et al., 2002). These
methods could be used to search approximate solutions of integral equations,
difference equations, integro-differential equations and partial
differential equations.

But the main difficulty in the application of perturbation methods is the
requirement of a small parameter or to install a small artificial parameter
in the equation. For this reason, the obtained solutions are restricted by a
validity range of physical parameters. Therefore, to overcome the
disadvantages come with the perturbation techniques, some methods have been
suggested by several authors (He, 2001; Mickens, 1987, 2005, 2006; Cooper,
2002; Hu and Xiong, 2003; He, 2012; Wang and He, 2008; Iqbal and Javed,
2011; Iqbal et al., 2010).

Parallel to these studies, recently a new perturbation-iteration algorithm
has been proposed by Aksoy, Pakdemirli and their co-workers, Aksoy and
Pakdemirli (2010); Pakdemirli et al. (2011); Aksoy et al. (2012). A previous
attempt of constructing root finding algorithms systematically, Pakdemirli
and Boyac\i (2007), Pakdemirli et al. (2007); Pakdemirli et al. (2008).
finally guided to generalization of the method to differential equations
also Aksoy and Pakdemirli, (2010); Pakdemirli et al. (2011); Aksoy et al.
(2012). In the new technique, an iterative algorithm is established on the
perturbation expansion. The method has been applied to first order equations
Pakdemirli et al. (2011) and Bratu type second order equations, Aksoy and
Pakdemirli, (2010) to obtain approximate solutions. Then the algorithms were
tested on some nonlinear heat equations also Aksoy et al. (2012). Finally,
the solutions of the Volterra and Fredholm type integral equations (Dolapci
et al., 2013) and ordinary differential equation and systems, \c{S}enol et
al. (2013) have given by the present method.

In this study, the previously developed new technique is applied to some
types of nonlinear fractional differential equations for the first time. To
obtain the approximate solutions of equations, the most basic
perturbation-iteration algorithm PIA(1,1) is employed by taking one
correction term in the perturbation expansion and correction terms of only
first derivatives in the Taylor series expansion, i.e. $n=1,m=1$.

Consider the following initial value problem.

\begin{equation}
F_{k}\left( D^{\alpha }u_{k},u_{j},\varepsilon ,t\right) =0  \label{4}
\end{equation}%
\[
k=1,2,\dots ,K
\]%
\[
j=1,2,\dots ,K
\]%
\begin{equation}
u^{(k)}(0)=c_{k},\qquad k=0,1,2,\ldots ,m-1,\quad m-1<\alpha \leq m
\label{5}
\end{equation}%
where $K$ is the number of the equation in the system and $D^{\alpha }$ is
the Caputo fractional derivative of order $\alpha $, which is defined by:%
\begin{equation}
D^{\alpha }u(t)=J^{m-\alpha }\left(\frac{d^{m}}{dt^{m}u(t)}\right),\quad m-1<\alpha
<m,\quad m\in \mathbf{N}  \label{6}
\end{equation}%
As more clearly the system could be expressed by:%
\begin{eqnarray}
{F_{1}} &{=}&F_{1}\left( u_{1}^{\left( \alpha \right) },u_{1},u_{2},{\dots ,u%
}_{k},\varepsilon ,t\right) =0  \nonumber \\
F_{2} &=&F_{2}\left( u_{2}^{\left( \alpha \right) },u_{1},u_{2},{\dots ,u}%
_{k},\varepsilon ,t\right) =0  \nonumber \\
&&\vdots  \nonumber \\
F_{K} &=&F_{K}\left( u_{k}^{\left( \alpha \right) },u_{1},u_{2},{\dots ,u}%
_{k},\varepsilon ,t\right) =0  \label{7}
\end{eqnarray}

In this method as $\varepsilon $ is the artificially introduced perturbation
parameter and we use only one correction term in the perturbation expansion.%
\begin{eqnarray}
u_{k,n+1} &=&u_{k,n}+\varepsilon (u_{c})_{k,n}  \nonumber \\
u_{k,n+1}^{\prime } &=&u_{k,n}^{\prime }+\varepsilon (u_{c}^{\prime })_{k,n}
\label{8}
\end{eqnarray}%
where subscript $n$ represents the $n.$iteration

Replacing $(\ref{8})$ into $(\ref{4})$ and writing in the Taylor Series
expansion for only first order derivatives in the neighborhood of $%
\varepsilon =0$ gives 
\begin{equation}
F_{K}=\sum_{m=0}^{M}{\frac{1}{m!}{\left[ {\left( \frac{d}{d\varepsilon }%
\right) }^{m}F_{K}\right] }_{\varepsilon =0}^{m}\varepsilon ^{m},\ \ \ \ }%
k=1,2,\dots ,K  \label{9}
\end{equation}%
for%
\begin{equation}
\frac{d}{d\varepsilon }=\frac{\partial u_{k,n+1}^{\left( \alpha \right) }}{%
\partial \varepsilon }\frac{\partial }{\partial u_{k,n+1}^{\left( \alpha
\right) }}+\sum_{j=1}^{K}{\left( \frac{\partial u_{j,n+1}}{\partial
\varepsilon }\frac{\partial }{\partial u_{j,n+1}}\right) }+\frac{\partial }{%
\partial \varepsilon }  \label{10}
\end{equation}

\bigskip This equation is defined for the $(n+1).$iteration equation%
\begin{equation}
F_{k}\left( u_{k,n+1}^{\left( \alpha \right) },u_{j,n+1},\varepsilon
,t\right) =0  \label{11}
\end{equation}

Replacing $(\ref{10})$ in $(\ref{9})$ yields our iteration equation:%
\begin{equation}
F_{K}=\sum_{m=0}^{M}\frac{1}{m!}{\left[ {\left( \frac{{\left(
(u_{c})_{k,n}\right) }^{\left( \alpha \right) }}{\partial \varepsilon }\frac{%
\partial }{\partial u_{k,n+1}^{\left( \alpha \right) }}%
+\sum_{j=1}^{K}(u_{c})_{j,n}{\frac{\partial }{\partial u_{j,n+1}}}+\frac{%
\partial }{\partial \varepsilon }\right) }^{m}F_{K}\right] }_{\varepsilon
=0}^{m}\varepsilon ^{m}=0,  \label{12}
\end{equation}%
\[
k=1,2,\dots ,K 
\]%
All derivatives are calculated at $\varepsilon =0.$

Beginning with an initial function $u_{0}$, first $(u_{c})_{k,n}$ 's has
been determined by the help of $(\ref{12})$. Then using Eq. $(\ref{8})$, $%
\left( n+1\right) .$ iteration solution could be found Iteration process is
repeated using $(\ref{12})$ and $(\ref{8})$ until achieving an acceptable
result.

\section{Application}

\begin{example}
Consider the following system of linear fractional differential equations Abdulaziz et al. 2008:
\end{example}%
\begin{eqnarray}
D^{\alpha }u_{1}\left( t\right) &=&u_{1}\left( t\right) +u_{2}\left( t\right)
\nonumber \\
D^{\beta }u_{1}\left( t\right) {\ } &=&-u_{1}\left( t\right) +u_{2}\left(
t\right)  \label{13}
\end{eqnarray}

\[
0<\alpha ,\beta \leq 1
\]%
given with the initial conditions $u_{1}\left( 0\right) =0$ and$\ \
u_{2}(0)=1$. The known exact solutions for $\alpha =\beta =1$ are%
\begin{equation}
u_{1}\left( t\right) =e^{t}sint  \label{14}
\end{equation}%
and 
\begin{equation}
u_{2}(t)=e^{t}cost  \label{15}
\end{equation}

Eq. $(\ref{13})$ can be rearranged in the following form with adding and
subtracting $u_{1,n}^{\prime }(t)$ and $u_{2,n}^{\prime }(t)$ to the
equation: 
\begin{equation}
\varepsilon \frac{d^{\alpha }u_{1}(t)}{{dt}^{\alpha }}+u_{1,n}^{\prime
}\left( t\right) -\varepsilon u_{1,n}^{\prime }\left( t\right) -\varepsilon
u_{1,n}\left( t\right) -\varepsilon u_{2,n}\left( t\right) =0  \label{16}
\end{equation}%
\begin{equation}
\varepsilon \frac{d^{\beta }u_{2}(t)}{{dt}^{\beta }}+u_{2,n}^{^{\prime
}}\left( t\right) -\varepsilon u_{2,n}^{^{\prime }}\left( t\right)
+\varepsilon u_{1,n}\left( t\right) -\varepsilon u_{2,n}\left( t\right) =0
\label{17}
\end{equation}%
where $\varepsilon $ is a small parameter. According to the iteration
formula for%
\begin{equation}
F\left( u_{1}^{^{\prime }},u_{1},\varepsilon \right) =\frac{1}{\Gamma
(1-\alpha )}\varepsilon \int_{0}^{t}{\frac{u_{1}^{\prime }(s)}{{(t-s)}%
^{\alpha }}ds+u_{1,n}^{^{\prime }}\left( t\right) -\varepsilon
u_{1,n}^{^{\prime }}\left( t\right) -\varepsilon u_{1,n}\left( t\right)
-\varepsilon u_{2,n}\left( t\right) }  \label{18}
\end{equation}%
\begin{equation}
F\left( u_{2}^{^{\prime }},u_{2},\varepsilon \right) =\frac{1}{\Gamma
(1-\beta )}\varepsilon \int_{0}^{t}{\frac{u_{2}^{\prime }(s)}{{(t-s)}^{\beta
}}ds+u_{2,n}^{^{\prime }}\left( t\right) -\varepsilon u_{2,n}^{^{\prime
}}\left( t\right) +\varepsilon u_{1,n}\left( t\right) -\varepsilon
u_{2,n}\left( t\right) }  \label{19}
\end{equation}%
terms in equation $(\ref{12})$ become%
\begin{eqnarray}
F &=&u_{1,n}^{^{\prime }}\left( t\right) ,\quad F_{u_{1}}=0,\quad
F_{u_{1}^{\prime }}=1,  \nonumber \\
F_{\varepsilon } &=&-u_{1,n}^{\prime }\left( t\right) -u_{1,n}\left(
t\right) -u_{2,n}\left( t\right) +\frac{1}{\Gamma (1-\alpha )}\int_{0}^{t}{%
\frac{u_{1,n}^{\prime }\left( s\right) }{{(t-s)}^{\alpha }}ds}  \label{20}
\end{eqnarray}%
for the iteration formula%
\begin{equation}
u_{1}^{^{\prime }}\left( t\right) +\frac{F_{u_{1}}}{F_{u_{1}^{\prime }}}%
u_{1}\left( t\right) =-\frac{F_{\varepsilon }+\frac{F}{\varepsilon }}{%
Fu_{1}^{\prime }}  \label{21}
\end{equation}%
and terms in equation $(\ref{12})$ 
\begin{eqnarray}
F &=&u_{2,n}^{\prime }\left( t\right) ,\quad F_{u_{2}}=0,\quad
F_{u_{2}^{\prime }}=1,  \nonumber \\
F_{\varepsilon } &=&-u_{2,n}^{^{\prime }}\left( t\right) +u_{1,n}\left(
t\right) -u_{2,n}\left( t\right) +\frac{1}{\Gamma (1-\beta )}\int_{0}^{t}{%
\frac{u_{2,n}^{\prime }\left( s\right) }{{(t-s)}^{\beta }}ds}  \label{22}
\end{eqnarray}%
for the iteration formula%
\begin{equation}
u_{2}^{^{\prime }}\left( t\right) +\frac{F_{u_{2}}}{F_{u_{2}^{^{\prime }}}}%
u_{2}\left( t\right) =-\frac{F_{\varepsilon }+\frac{F}{\varepsilon }}{%
F_{u_{2}^{\prime }}}  \label{23}
\end{equation}%
Notice that inserting the small parameter $\varepsilon $ as a coefficient of
the integral term simplifies the system and makes it easy to solve.

After writing these terms in the iteration formulas, the system gives the
following differential equations. 
\begin{equation}
u_{1,n}\left( t\right) +u_{2,n}\left( t\right) +\frac{(-1+\varepsilon )u_{1},%
{_{n}}^{^{\prime }}(t)}{\varepsilon }=\frac{1}{\Gamma \left( 1-\alpha
\right) }\int_{0}^{t}{\frac{u_{1,n}^{^{\prime }}\left( s\right) }{{\left(
t-s\right) }^{\alpha }}ds+}({u_{c})_{1,n}^{\prime }(t)}  \label{24}
\end{equation}%
and 
\begin{equation}
{-}u_{1,n}\left( t\right) +u_{2,n}\left( t\right) +\frac{(-1+\varepsilon
)u_{2,n}^{^{\prime }}(t)}{\varepsilon }=\frac{1}{\Gamma \left( 1-\beta
\right) }\int_{0}^{t}{\frac{u_{2,n}^{^{\prime }}\left( s\right) }{{\left(
t-s\right) }^{\beta }}ds+}({u_{c})_{2,n}^{\prime }(t)}  \label{25}
\end{equation}

When employing the iteration formula $(\ref{12})$, we start with an initial
function compatible to the boundary condition and we determine coefficients
from the boundary condition at each step. Beginning with the initial
functions 
\[
u_{1,0}\left( t\right) =0~ \text{ and }~ u_{2,0}\left( t\right) =1
\]%
and using the iteration formula, we get the following successive approximate
solutions at each step for $n=0,1,2,\dots $%
\begin{eqnarray}
u_{1,1}\left( t\right)  &=&t  \nonumber \\
u_{2,1}\left( t\right)  &=&1+t  \nonumber \\
u_{1,2}\left( t\right)  &=&t\left( 2+t-\frac{t^{1-\alpha }}{\Gamma (3-\alpha
)}\right)   \nonumber \\
u_{2,2}\left( t\right)  &=&1+2t-\frac{t^{2-\beta }}{\Gamma (3-\beta )} 
\nonumber \\
u_{1,3}\left( t\right)  &=&\frac{1}{3}t\left( 9+t(9+t)+\frac{3t^{2-2\alpha }%
}{\Gamma \left( 4-2\alpha \right) }-\frac{9t^{1-\alpha }(3+t-\alpha )}{%
\Gamma \left( 4-\alpha \right) }-\frac{3t^{2-\beta }}{\Gamma \left( 4-\beta
\right) }\right)   \nonumber \\
u_{2,3}\left( t\right)  &=&1+3t-\frac{t^{3}}{3}+\frac{t^{3-\alpha }}{\Gamma
\left( 4-\alpha \right) }+t^{2-2\beta }\left( \frac{t}{\Gamma \left(
4-2\beta \right) }-\frac{t^{\beta }\left( 9+t-3\beta \right) }{\Gamma \left(
4-\beta \right) }\right)   \label{26}
\end{eqnarray}%
and so on. Following in this manner the fifth iteration solutions $u_{1,5}(t)
$ and $u_{2,5}(t)$ were calculated. Due to brevity reasons, higher
iterations are not given here. It is easy to calculate other iterations up
to any order by the help of a symbolic calculation software such as
Mathematica . In Figure 1. and 2. we compare our $u_{1,5}(t)$ and $u_{2,5}(t)
$ solutions for different values of $\alpha $ and $\beta $ and in Figure 3.
and 4. with the exact solution graphically. In Table 1. and 2. some of our
iterations are compared with the exact solution. The results show that the
proposed method can give successful approximations.

{\scriptsize 
\begin{table}[tbp]
\caption{Numerical results of Example 1. for some values of $u_{1}(t)$.}
\begin{center}
{\scriptsize 
\begin{tabular}{cccccccc}
\hline
\multicolumn{8}{c}{$\alpha =\beta =1$} \\ \hline
$t$ & $u_{1,1}(t)$ & $u_{1,2}(t)$ & $u_{1,3}(t)$ & $u_{1,4}(t)$ & $u_{1,5}(t)$ & $Exact$ $Solution$ & $Absolute$ $Error$ \\ \hline
\textbf{0.0} & 0.000000 & 0.000000 & 0.000000 & 0.000000 & 0.000000 & 
0.000000 & 0.000000 \\ 
\textbf{0.1} & 0.100000 & 0.110000 & 0.110333 & 0.110333 & 0.110332 & 
0.110332 & 1.12698E-8 \\ 
\textbf{0.2} & 0.200000 & 0.240000 & 0.242666 & 0.242666 & 0.242656 & 
0.242655 & 7.31405E-7 \\ 
\textbf{0.3} & 0.300000 & 0.390000 & 0.399000 & 0.399000 & 0.398919 & 
0.398910 & 8.44622E-6 \\ 
\textbf{0.4} & 0.400000 & 0.559999 & 0.581333 & 0.581333 & 0.580991 & 
0.580943 & 0.00004809 \\ 
\textbf{0.5} & 0.500000 & 0.750000 & 0.791666 & 0.791666 & 0.790625 & 
0.790439 & 0.00018591 \\ 
\textbf{0.6} & 0.600000 & 0.960000 & 1.032000 & 1.031999 & 1.029408 & 
1.028845 & 0.00056233 \\ 
\textbf{0.7} & 0.700000 & 0.190000 & 1.304333 & 1.304333 & 1.298731 & 
1.297295 & 0.00143589 \\ 
\textbf{0.8} & 0.800000 & 1.440000 & 1.610666 & 1.610666 & 1.599744 & 
1.596505 & 0.00323866 \\ 
\textbf{0.9} & 0.900000 & 1.710000 & 1.953000 & 1.952999 & 1.933316 & 
1.926673 & 0.00664370 \\ 
\textbf{1.0} & 1.000000 & 2.000000 & 2.333333 & 2.333333 & 2.300000 & 
2.287355 & 0.01264471 \\ \hline
\end{tabular}
}
\end{center}
\end{table}
}

{\scriptsize 
\begin{table}[tbp]
\caption{Numerical results of Example 1. for some values of $u_{2}(t)$.}
\begin{center}
{\scriptsize 
\begin{tabular}{cccccccc}
\hline
\multicolumn{8}{c}{$\alpha =\beta =1$} \\ \hline
$t$ & $u_{2,1}(t)$ & $u_{2,2}(t)$ & $u_{2,3}(t)$ & $u_{2,4}(t)$ & $u_{2,5}(t)$ & $Exact$ $Solution$ & $Absolute$ $Error$ \\ \hline
\textbf{0.0} & 1.000000 & 1.000000 & 1.000000 & 1.000000 & 1.000000 & 
1.000000 & 0.00000000 \\ \hline
\textbf{0.1} & 1.100000 & 1.099999 & 1.099666 & 1.099650 & 1.099649 & 
1.099649 & 1.6274E-10 \\ \hline
\textbf{0.2} & 1.200000 & 1.200000 & 1.197333 & 1.197066 & 1.197055 & 
1.197056 & 2.13559E-8 \\ \hline
\textbf{0.3} & 1.300000 & 1.300000 & 1.291000 & 1.289650 & 1.289568 & 
1.289569 & 3.74045E-7 \\ \hline
\textbf{0.4} & 1.400000 & 1.400000 & 1.378666 & 1.374399 & 1.374058 & 
1.374061 & 2.87222E-6 \\ \hline
\textbf{0.5} & 1.500000 & 1.500000 & 1.458333 & 1.447916 & 1.446874 & 
1.446889 & 0.00001403 \\ \hline
\textbf{0.6} & 1.600000 & 1.600000 & 1.528000 & 1.506399 & 1.503807 & 
1.503859 & 0.00005154 \\ \hline
\textbf{0.7} & 1.700000 & 1.700000 & 1.585666 & 1.545650 & 1.540047 & 
1.540803 & 0.00015535 \\ \hline
\textbf{0.8} & 1.800000 & 1.800000 & 1.629333 & 1.561066 & 1.550144 & 
1.550549 & 0.00040529 \\ \hline
\textbf{0.9} & 1.900000 & 1.900000 & 1.659999 & 1.547650 & 1.527966 & 
1.528913 & 0.00094681 \\ \hline
\textbf{1.0} & 2.000000 & 2.000000 & 1.666666 & 1.500000 & 1.466666 & 
1.468693 & 0.00202727 \\ \hline
\end{tabular}
}
\end{center}
\end{table}
}

\begin{figure}[h]
$\centering\includegraphics[width=0.75\textwidth]{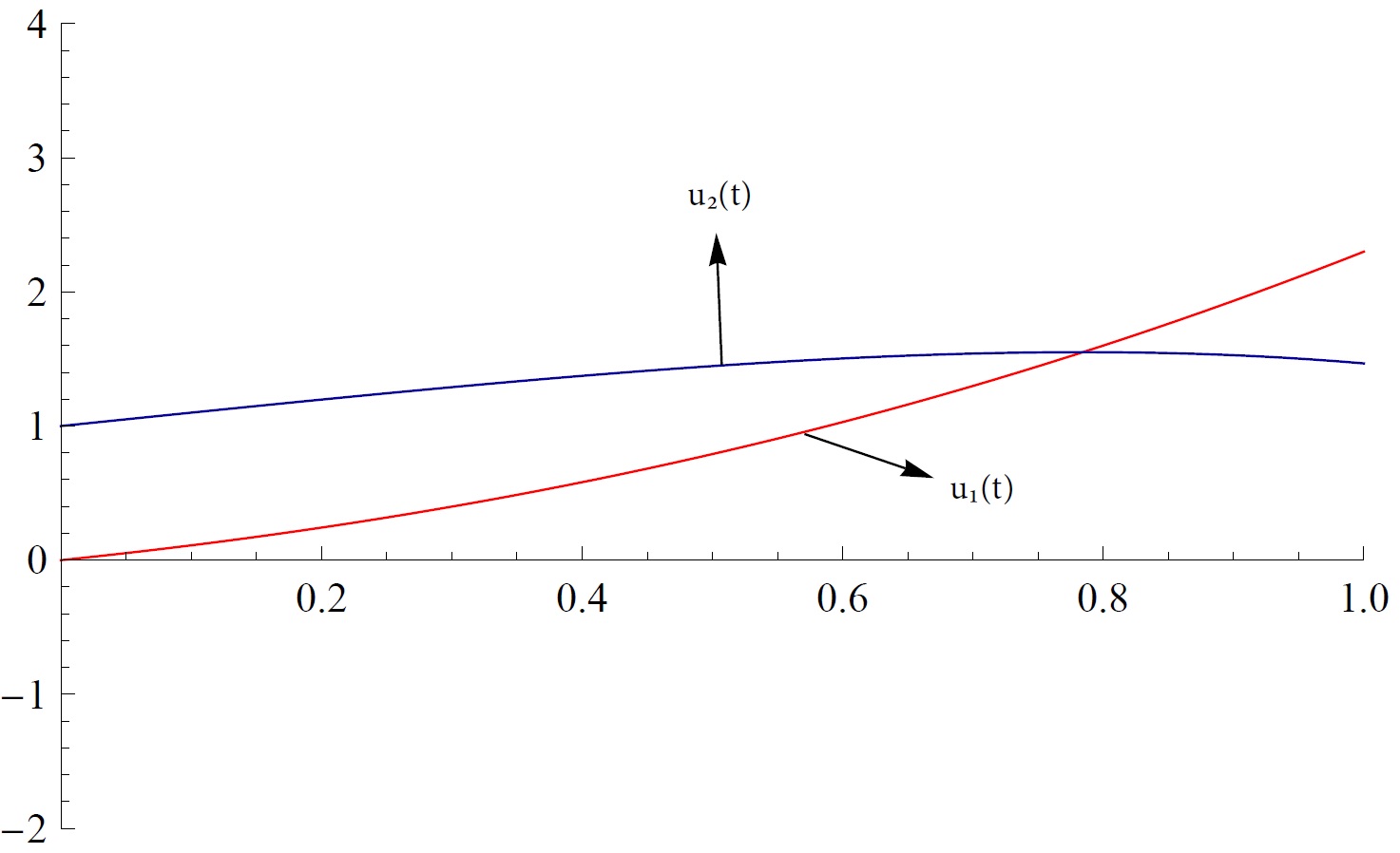}$%
\caption{Comprasion of the $u_{1,5}(t)$ and $u_{2,5}(t)$ of Example 1. for $%
\protect\alpha =\protect\beta =1$.}
\label{fig:Figure1}
\end{figure}

\begin{figure}[h]
\centering
\includegraphics[width=0.75\textwidth]{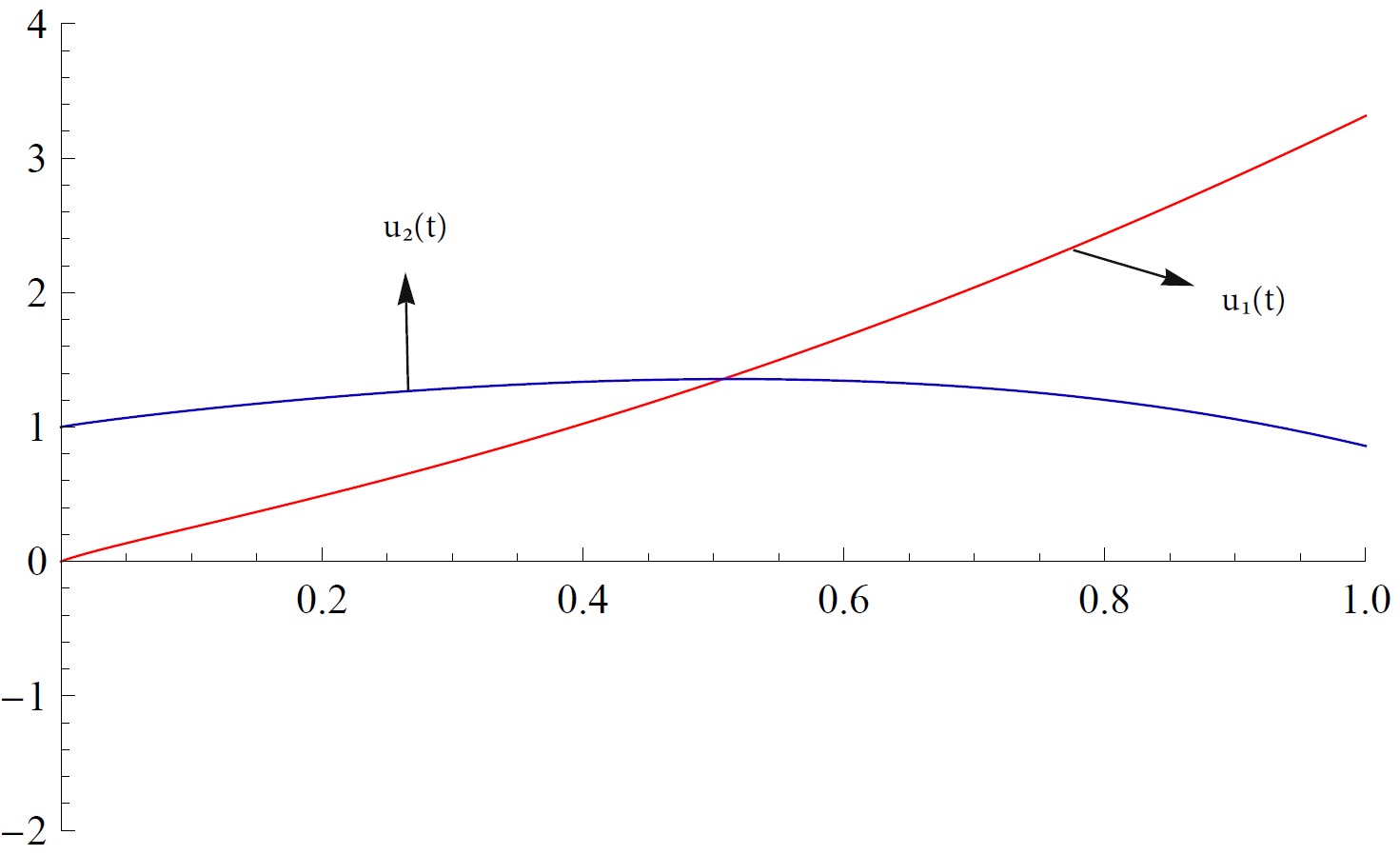}
\caption{Comprasion of the $u_{1,5}(t)$ and $u_{2,5}(t)$ of Example 1. for $%
\protect\alpha =0.7$ and $\protect\beta =0.9$.}
\label{fig:Figure2}
\end{figure}

\begin{figure}[h]
\centering
\includegraphics[width=0.75\textwidth]{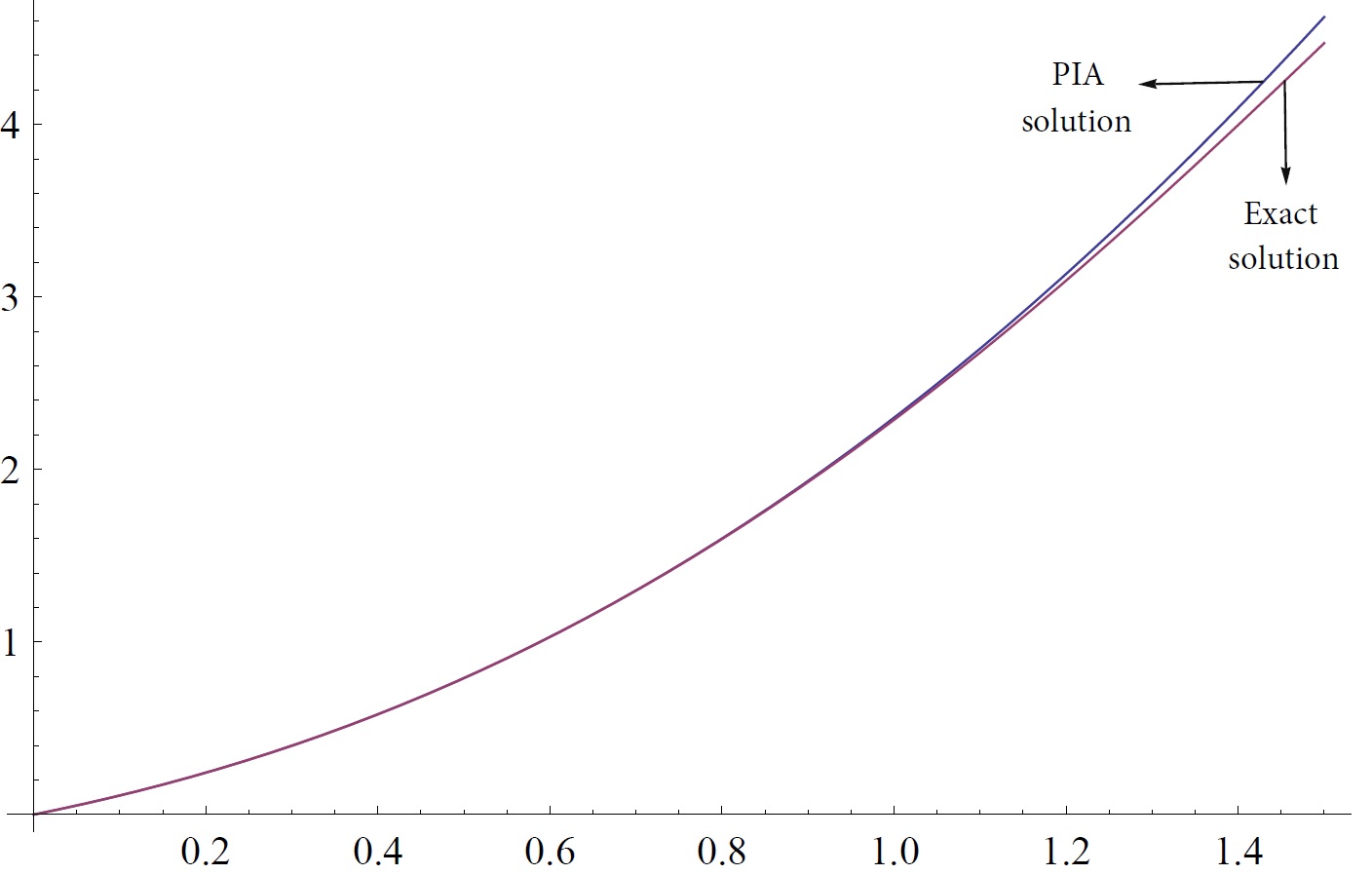}
\caption{Comprasion of the exact solution and the PIA solution $u_{1,5}(t)$
of Example 1. for $\protect\alpha =\protect\beta =1$.}
\label{fig:Figure3}
\end{figure}

\begin{figure}[h]
\centering
\includegraphics[width=0.75\textwidth]{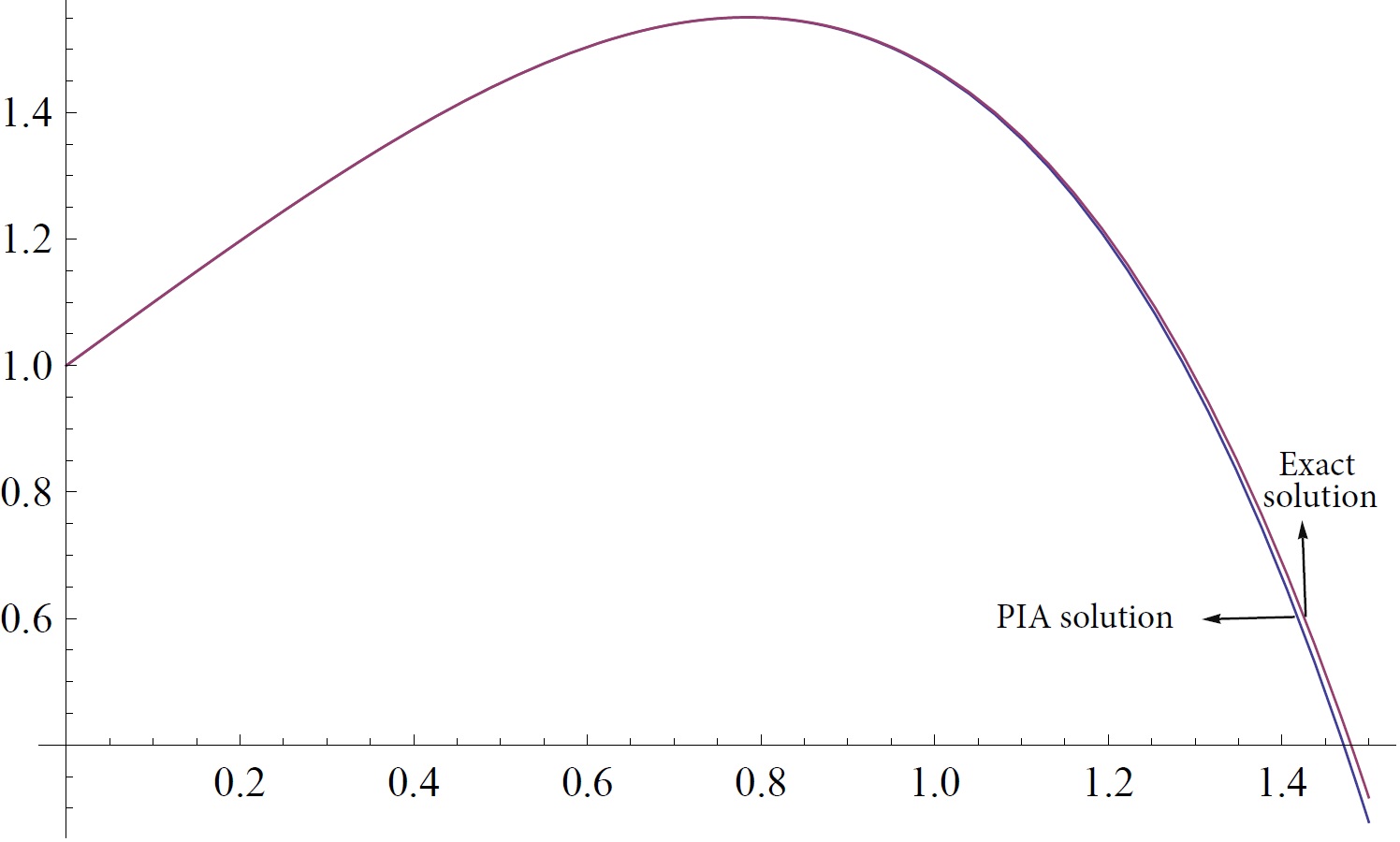}
\caption{Comprasion of the exact solution and the PIA solution $u_{2,5}(t)$
of Example 1. for $\protect\alpha =\protect\beta =1$.}
\label{fig:Figure4}
\end{figure}

\newpage

\begin{example}
For the second example consider the following system of nonlinear fractional differential equations Zurigat et al. 2010: 
\end{example}%
\begin{eqnarray}
D^{\alpha _{1}}{u}_{1}{\left( t\right) \ } &=&\frac{1}{2}{u}_{1}\left(
t\right)  \nonumber \\
D^{\alpha _{2}}u_{2}{\left( t\right) \ } &=&u_{2}\left( t\right) +{u}_{1}^{2}%
{\left( t\right) }  \label{27-1}
\end{eqnarray}

\[
0<\alpha _{1},\alpha _{2}\leq 1
\]%
given with the initial conditions ${u}_{1}\left( 0\right) =1$ and$\
u_{2}(0)=0$. The exact solutions, when $\alpha _{1}=\alpha _{2}=1$, are

\begin{equation}
{u}_{1}\left( t\right) =e^{\frac{t}{2}}  \label{27}
\end{equation}%
and

\begin{equation}
u_{2}(t)={te}^{t}  \label{28}
\end{equation}

In the system, if add and subtract $u_{1,n}^{\prime }(t)$ and $%
u_{2,n}^{\prime }(t)$ respectively, the system could be rewritten in the
following form: 
\begin{equation}
\varepsilon D^{\alpha _{1}}{u}_{1}{\left( t\right) \ \ }+{u}_{1,n}^{\prime
}\left( t\right) -\varepsilon {u}_{1,n}^{\prime }\left( t\right) -\frac{1}{2}%
\varepsilon {u}_{1,n}\left( t\right) =0  \label{29}
\end{equation}%
\begin{equation}
\varepsilon D^{\alpha _{2}}u_{2}{\left( t\right) \ }+u_{2,n}^{\prime }\left(
t\right) -\varepsilon u_{2,n}^{\prime }\left( t\right) -\varepsilon
u_{2,n}\left( t\right) -\varepsilon {u}_{1,n}^{2}\left( t\right) {\ }=0
\label{30}
\end{equation}%
where $\varepsilon $ is a small parameter. For%
\begin{equation}
F\left( {u}_{1}^{^{\prime }},{u}_{1},\varepsilon \right) =\frac{1}{\Gamma
(1-\alpha _{1})}\varepsilon \int_{0}^{t}{\frac{{u}_{1,n}^{\prime }(s)}{{(t-s)%
}^{\alpha _{1}}}ds+{u}_{1,n}^{\prime }(t)-\varepsilon {u}_{1,n}^{\prime }(t)-%
\frac{1}{2}\varepsilon {u}_{1,n}\left( t\right) }  \label{31}
\end{equation}%
\begin{equation}
F\left( u_{2}^{^{\prime }},u_{2},\varepsilon \right) =\frac{1}{\Gamma
(1-\alpha _{2})}\varepsilon \int_{0}^{t}{\frac{u_{2,n}^{\prime }(s)}{{(t-s)}%
^{\alpha _{2}}}ds+{u}_{2,n}^{\prime }(t)-\varepsilon {u}_{2,n}^{\prime
}(t)-\varepsilon u_{2,n}\left( t\right) -\varepsilon {u}_{1,n}^{2}{\left(
t\right) }}  \label{32}
\end{equation}%
and terms in equation $(\ref{12})$ become 
\begin{eqnarray}
F &=&{u}_{1,n}^{\prime }(t),\quad F_{{u}_{1}}=0,\quad F_{{u}_{1}^{\prime
}}=1,  \nonumber \\
F_{\varepsilon } &=&-{u}_{1,n}^{\prime }(t)-\frac{{u}_{1,n}(t)}{2}+\frac{1}{%
\Gamma (1-\alpha _{1})}\int_{0}^{t}{\frac{{u}_{1,n}^{\prime }(t)}{{(t-s)}%
^{\alpha _{1}}}ds}  \label{33}
\end{eqnarray}%
and%
\begin{eqnarray}
F &=&u_{2,n}^{\prime }(t),\quad F_{u_{2}}=0,\quad F_{u_{2}^{\prime }}=1, 
\nonumber \\
F_{\varepsilon } &=&-u_{2,n}^{\prime }(t)-u_{2,n}(t)-{u}_{1,}^{2}{_{n}}(t)+%
\frac{1}{\Gamma (1-\alpha _{2})}\int_{0}^{t}{\frac{u_{2,n}^{\prime }(s)}{{%
(t-s)}^{\alpha _{2}}}ds}  \label{34}
\end{eqnarray}%
\qquad \qquad\ \ After writing these terms in the iteration formula we
obtain the following differential equations: 
\begin{equation}
2\left( \frac{\int_{0}^{t}{{(-s+t)}^{-\alpha _{1}}{u}_{1,n}^{\prime }(s)ds}}{%
\Gamma (1-\alpha _{1})}+({u_{c})_{1,n}^{\prime }(t)}-\frac{(-1+\varepsilon ){%
({u}_{1,n})}^{^{\prime }}(t)}{\varepsilon }\right) ={u}_{1,n}\left( t\right) 
\label{35}
\end{equation}%
and 
\begin{equation}
u_{2,n}\left( t\right) +u_{1,n}^{2}{\left( t\right) }+\frac{\left(
-1+\varepsilon \right) {\left( u_{2,n}\right) }^{^{\prime }}\left( t\right) 
}{\varepsilon }=\frac{\int_{0}^{t}{{(-s+t)}^{-\alpha _{2}}{(u_{2,n})}%
^{^{\prime }}s)ds}}{\Gamma (1-\alpha _{2})}+({u_{c})_{2,n}^{\prime }(t)}
\label{36}
\end{equation}

Beginning with the initial functions 
\[
{u}_{1,0}\left( t\right) =1~\text{and}~u_{2,0}\left( t\right) =0 
\]%
and using the iteration formula, the following successive approximate
solutions are obtained for $n=0,1,2,\dots $

\begin{eqnarray}
u_{1,1}(t) &=&1+\frac{t}{2}  \nonumber \\
u_{2,1}(t) &=&t  \nonumber \\
u_{1,2}(t) &=&\frac{1}{8}\left( 8+8t+t^{2}+\frac{4t^{2-\alpha _{1}}}{\Gamma
(2-\alpha _{1})(-2+\alpha _{1})}\right)  \nonumber \\
u_{2,2}\left( t\right) &=&2t+t^{2}+\frac{t^{3}}{12}+\frac{t^{2-\alpha _{2}}}{%
\Gamma (2-\alpha _{2})\left( -2+\alpha _{2}\right) }  \nonumber \\
u_{1,3}\left( t\right) &=&\frac{1}{48}\left( 48+72t+18t^{2}+t^{3}+\frac{%
24t^{3-2\alpha _{1}}}{\Gamma (4-2\alpha _{1})}-\frac{24t^{2-\alpha
_{1}}\left( 9+t-3\alpha _{1}\right) }{\Gamma (4-\alpha _{1})}\right) 
\nonumber \\
u_{2,3}\left( t\right) &=&3t+3t^{2}+\frac{5t^{3}}{6}+\frac{t^{4}}{12}+\frac{%
t^{5}}{320}+\frac{t^{3-2\alpha _{2}}}{\Gamma (4-2\alpha _{2})}-\frac{%
t^{5-2\alpha _{1}}}{4\Gamma {(3-\alpha _{1})}^{2}\left( -5+2\alpha
_{1}\right) }  \nonumber \\
&&-\frac{t^{3-\alpha _{1}}\left( 4\left( 40+3t\left( 10+t\right) \right)
+\alpha _{1}\left( -72-t\left( 64+7t\right) +\left( 8+t\left( 8+t\right)
\right) \alpha _{1}\right) \right) }{8\Gamma \left( 6-\alpha _{1}\right) } 
\nonumber \\
&&-\frac{t^{2-\alpha _{2}}(72+t(24+t)+6\alpha _{2}(-7-t+\alpha _{2}))}{%
2\Gamma (5-\alpha _{2})}  \label{37}
\end{eqnarray}

and so on. In the same manner the fourth iteration solutions $u_{1,4}(t)$
and $u_{2,4}(t)$ are calculated. Again we compared our results In Figures
5., 6., 7. and 8. and in Table 3. and 4. with the exact solutions.

{\scriptsize 
\begin{table}[tbp]
\caption{Numerical results of Example 2. for some values of $u_{1}(t)$.}
\begin{center}
{\scriptsize 
\begin{tabular}{ccccccc}
\hline
\multicolumn{7}{c}{$\alpha _{1} =\alpha _{2} =1$} \\ \hline
$t$ & $u_{1,1}(t)$ & $u_{1,2}(t)$ & $u_{1,3}(t)$ & $u_{1,4}(t)$ & $Exact$ $Solution$ & $Absolute$ $Error$ \\ \hline
\textbf{0.0} & 1.000000 & 1.000000 & 1.000000 & 1.000000 & 1.000000 & 
0.0000000 \\ 
\textbf{0.1} & 1.050000 & 1.051250 & 1.051270 & 1.051271 & 1.051271 & 
2.6260E-9 \\ 
\textbf{0.2} & 1.100000 & 1.105000 & 1.105166 & 1.105170 & 1.105170 & 
8.4742E-8 \\ 
\textbf{0.3} & 1.150000 & 1.161250 & 1.161812 & 1.161833 & 1.161834 & 
6.4897E-7 \\ 
\textbf{0.4} & 1.200000 & 1.220000 & 1.221333 & 1.221400 & 1.221402 & 
2.7581E-6 \\ 
\textbf{0.5} & 1.250000 & 1.281250 & 1.283854 & 1.284016 & 1.284025 & 
8.4896E-6 \\ 
\textbf{0.6} & 1.300000 & 1.345000 & 1.349500 & 1.349837 & 1.349858 & 
0.0000213 \\ 
\textbf{0.7} & 1.350000 & 1.411250 & 1.418395 & 1.419021 & 1.419067 & 
0.0000464 \\ 
\textbf{0.8} & 1.400000 & 1.480000 & 1.490666 & 1.491733 & 1.491824 & 
0.0000913 \\ 
\textbf{0.9} & 1.450000 & 1.551250 & 1.566437 & 1.568146 & 1.568312 & 
0.0001660 \\ 
\textbf{1.0} & 1.500000 & 1.625000 & 1.645833 & 1.648437 & 1.648721 & 
0.0002837 \\ \hline
\end{tabular}
}
\end{center}
\end{table}
}

{\scriptsize 
\begin{table}[tbp]
\caption{Numerical results of Example 2. for some values of $u_{2}(t)$.}
\begin{center}
{\scriptsize 
\begin{tabular}{ccccccc}
\hline
\multicolumn{7}{c}{$\alpha _{1} =\alpha _{2} =1$} \\ \hline
$t$ & $u_{2,1}(t)$ & $u_{2,2}(t)$ & $u_{2,3}(t)$ & $u_{2,4}(t)$ & $Exact$ $Solution$ & $Absolute$ $Error$ \\ \hline
\textbf{0.0} & 0.000000 & 0.000000 & 0.000000 & 0.000000 & 0.000000 & 
0.0000000 \\ \hline
\textbf{0.1} & 0.100000 & 0.110083 & 0.110505 & 0.110516 & 0.110517 & 
2.4666E-7 \\ \hline
\textbf{0.2} & 0.200000 & 0.240666 & 0.244084 & 0.244272 & 0.244280 & 
8.12862-6 \\ \hline
\textbf{0.3} & 0.300000 & 0.392250 & 0.403929 & 0.404894 & 0.404957 & 
0.0000635 \\ \hline
\textbf{0.4} & 0.400000 & 0.565333 & 0.593365 & 0.596453 & 0.596729 & 
0.0002760 \\ \hline
\textbf{0.5} & 0.500000 & 0.760416 & 0.815852 & 0.823492 & 0.824360 & 
0.0008683 \\ \hline
\textbf{0.6} & 0.600000 & 0.978000 & 1.074993 & 1.091043 & 1.0932712 & 
0.0022277 \\ \hline
\textbf{0.7} & 0.700000 & 1.218583 & 1.374530 & 1.404661 & 1.409626 & 
0.0049654 \\ \hline
\textbf{0.8} & 0.800000 & 1.482666 & 1.718357 & 1.770446 & 1.780432 & 
0.0099863 \\ \hline
\textbf{0.9} & 0.900000 & 1.770750 & 2.110517 & 2.195074 & 2.213642 & 
0.0185684 \\ \hline
\textbf{1.0} & 1.000000 & 2.083333 & 2.555208 & 2.685825 & 2.718281 & 
0.0324559 \\ \hline
\end{tabular}
}
\end{center}
\end{table}
}

{\scriptsize 
\begin{figure}[h]
{\scriptsize \centering
\includegraphics[width=0.75\textwidth]{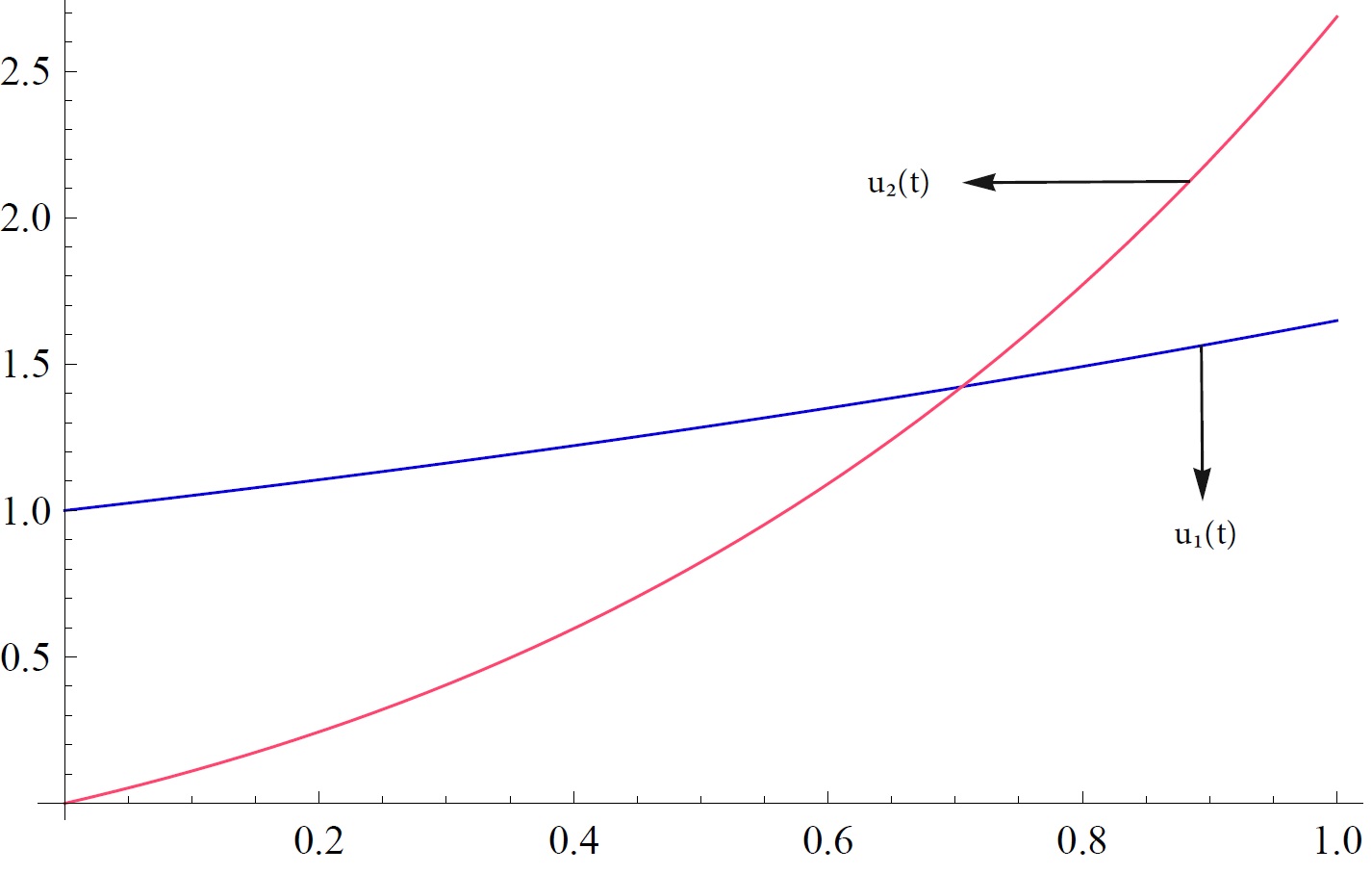}  }
\caption{Comprasion of the $u_{1,4}(t)$ and $u_{2,4}(t)$ of Example 2. for $%
\protect\alpha _{1}=\protect\alpha _{2}=1$.}
\label{fig:Figure5}
\end{figure}
}

{\scriptsize 
\begin{figure}[h]
{\scriptsize \centering
\includegraphics[width=0.75\textwidth]{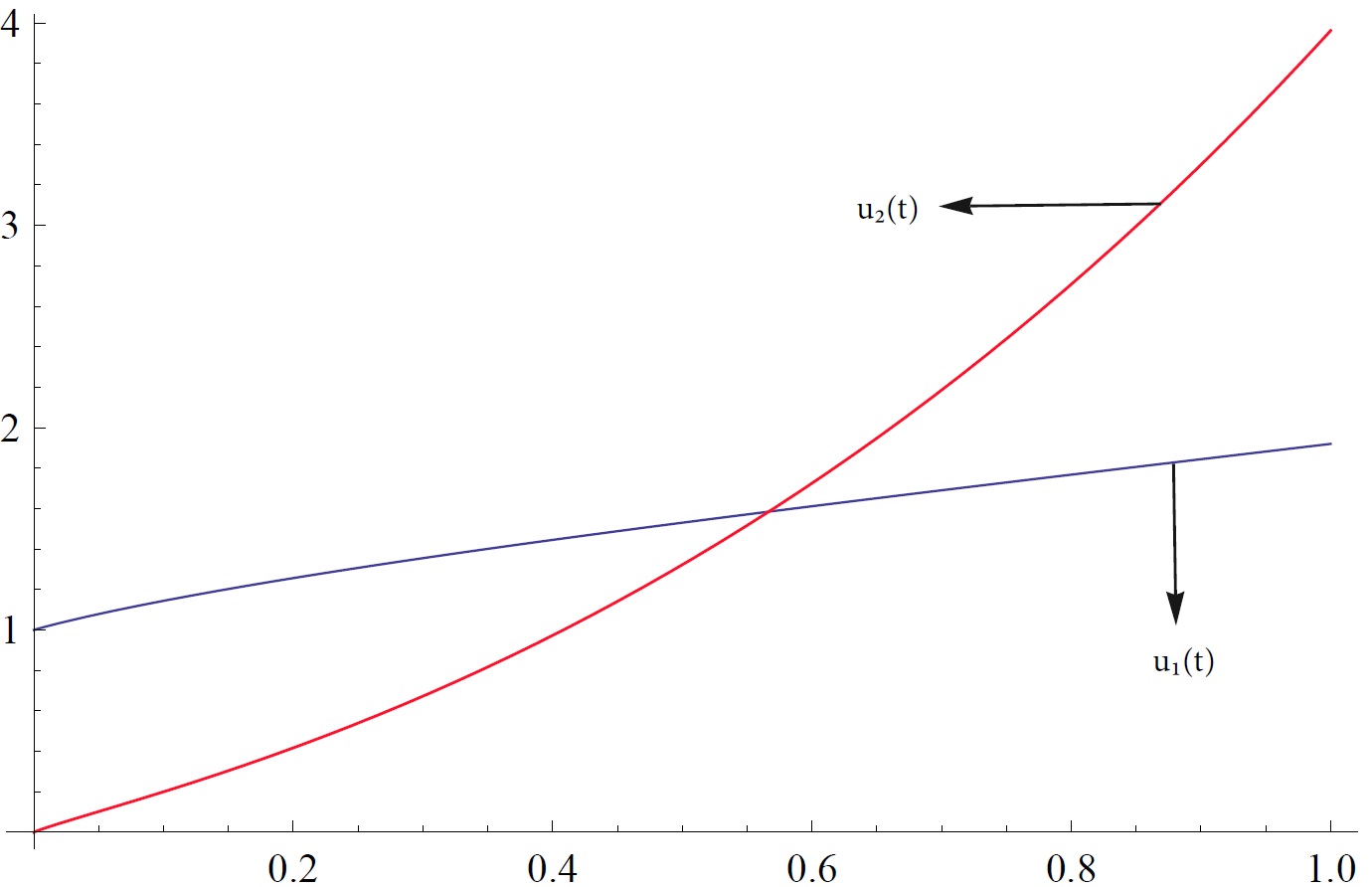}  }
\caption{Comprasion of the $u_{1,4}(t)$ and $u_{2,4}(t)$ of Example 2. for $%
\protect\alpha _{1}=0.5$ and $\protect\alpha _{2}=0.8$.}
\label{fig:Figure6}
\end{figure}
}

{\scriptsize 
\begin{figure}[h]
{\scriptsize \centering
\includegraphics[width=0.75\textwidth]{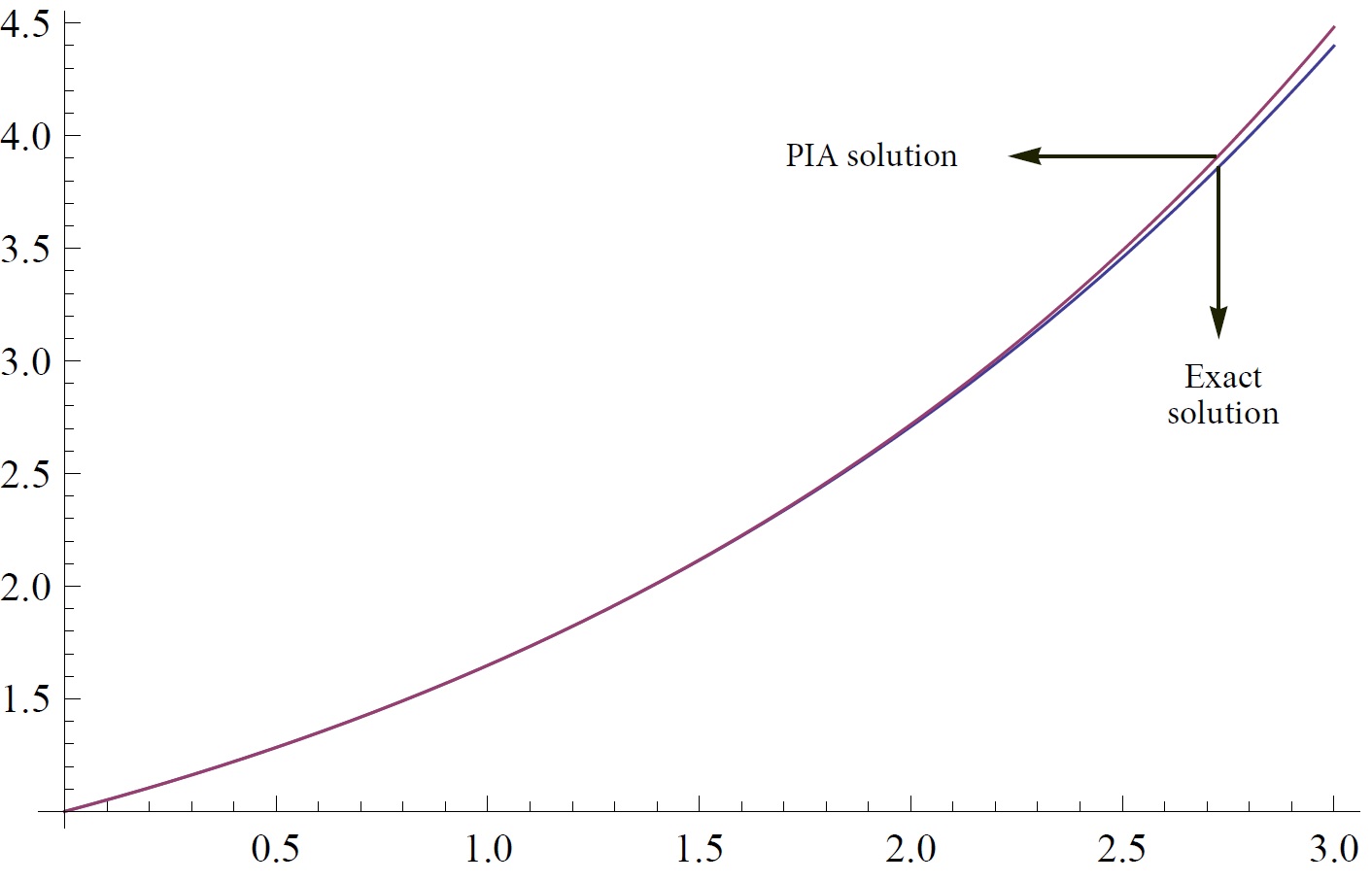}  }
\caption{Comprasion of the exact solution and the PIA solution $u_{1,4}(t)$
of Example 2. for $\protect\alpha _{1}=\protect\alpha _{2}=1$.}
\label{fig:Figure7}
\end{figure}
}

{\scriptsize 
\begin{figure}[h]
{\scriptsize \centering
\includegraphics[width=0.75\textwidth]{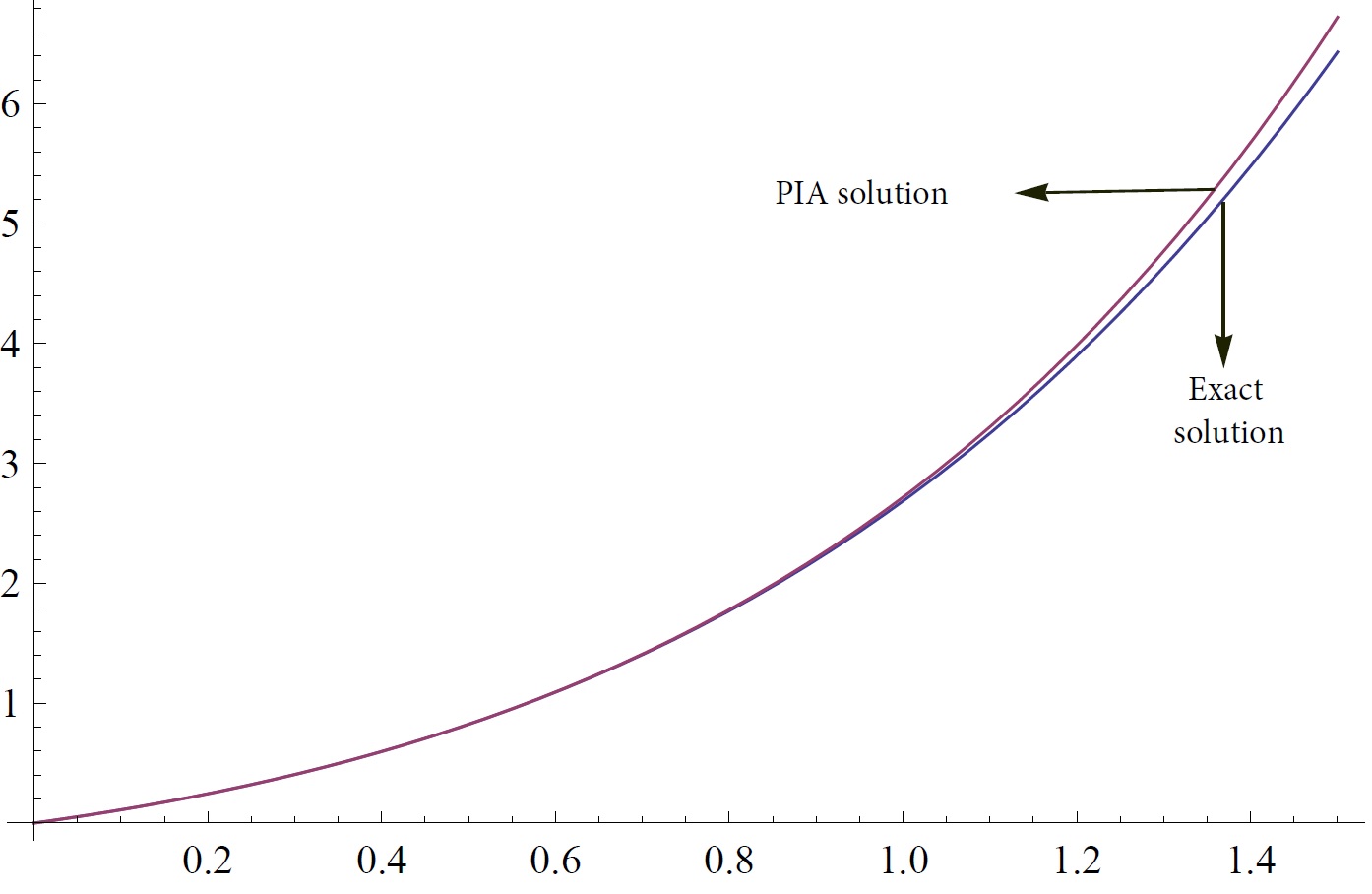}  }
\caption{Comprasion of the exact solution and the PIA solution $u_{2,4}(t)$
of Example 2. for $\protect\alpha _{1}=\protect\alpha _{2}=1$.}
\label{fig:Figure8}
\end{figure}
}

\section{Conclusion}

In this paper we have applied previously developed numerical method
so-called Perturbation-Iteration Algorithm (PIA) to find approximate
solutions of some system of nonlinear Factional Differential Equations for
the first time. The numerical results obtained in this study show that
method PIA is a remarkably successful numerical technique for solving system
of FDEs. We expect that the present method could used to calculate the
approximate solutions of other types of fractional differential equations
such as fractional integro-differential equations and fractional partial
differential equations. Our next study will be about these types of
equations.

\section{Authors’ contributions}

All authors contributed in development of the manuscript and solving problems. All authors read and approved the final manuscript.

\section{Competing interests}

The authors declare that they have no competing interest.

\end{document}